\documentclass{elsart3-1}

 \usepackage{amssymb}

\usepackage[english,francais]{babel}

\newtheorem{theorem}{Theorem}[section]

\newtheorem{e-proposition}[theorem]{Proposition}

\newtheorem{e-definition}[theorem]{Definition\rm}

\newtheorem{theoreme}{Th\'eor\`eme}[section]

\newtheorem{definition}[theoreme]{D\'efinition\rm}
\newtheorem{remarque}[theoreme]{\it Remarque}
\newtheorem{exemple}{\it Exemple\/}
\renewcommand{\theequation}{\arabic{equation}}
\def\om{{\omega}}

\def\si{{\sigma}}
\def\Si{{\Sigma}}

\def\ep{{\varepsilon}}

\def\phi{{\varphi}}

\DeclareMathAlphabet{\doba}{U}{msb}{m}{n}
\gdef\mC{\doba{C}}

\gdef\mR{\doba{R}}
\gdef\mS{\doba{S}}

\def\lamin{\lambda_{\rm min}^+}
\def\laminm{\lambda_{\rm min}^-}

\def\Vol{{\mathop{\rm Vol}}}
\def\Scal{{\mathop{\rm Scal}}}
\def\Ric{{\mathop{\rm Ric}}}

\def\Hom{{\mathop{\rm Hom}}}

\let\<\langle
\let\>\rangle

\long\def\komment#1{}

  \def\og{\leavevmode\raise.3ex\hbox{$\scriptscriptstyle\langle\!\langle$~}}
\def\fg{\leavevmode\raise.3ex\hbox{~$\!\scriptscriptstyle\,\rangle\!\rangle$}}

\begin{document}

\begin{frontmatter}

\selectlanguage{francais}
\title{Un probl\`eme de type Yamabe sur les vari\'et\'es spinorielles
  compactes }

\vspace{-2.6cm}
\selectlanguage{english}
\title{ A Yamabe type problem on compact spin manifolds}

\author[bernd]{Bernd Ammann}
\author[emmanuel]{Emmanuel Humbert}
\author{Bertrand Morel$^\mathrm{b}$}

\address[bernd]{Courriel: bernd.ammann@gmx.de}
\address[emmanuel]{Institut \'Elie Cartan BP 239
Universit\'e de Nancy 1
54506 Vandoeuvre-l\`es-Nancy Cedex
France   Courriel: humbert@iecn.u-nancy.fr}

\maketitle
\begin{abstract}
Let $(M,g,\si)$ be a compact spin manifold of dimension $n \geq
2$.
Let $\lambda_1^+(\tilde{g})$ be the smallest
positive eigenvalue of the Dirac operator in the metric
$\tilde{g} \in [g]$ conformal to  $g$. We then define
$\lamin(M,[g],\si) = \inf_{\tilde{g} \in [g] } \lambda_1^+(\tilde{g})
\Vol(M,\tilde{g})^{1/n} $.  We show that
$0< \lamin(M,[g],\si) \leq \lamin(\mS^n)$.
We find sufficient conditions for which we obtain strict inequality
$\lamin(M,[g],\si) < \lamin(\mS^n)$.
This strict inequality has applications to conformal spin geometry.\\[1mm]

\vskip 0.5\baselineskip

\selectlanguage{francais}
\noindent{\bf R\'esum\'e}
\vskip 0.5\baselineskip
\noindent
  Soit $(M,g,\si)$ une vari\'et\'e spinorielle compacte de dimension $n \geq
2$.
On note
$\lambda_1^+(\tilde{g})$  la plus petite
valeur propre $>0$ de
l'op\'erateur de Dirac dans la m\'etrique $\tilde{g} \in [g]$ conforme \`a
$g$.  On d\'efinit
$\lamin(M,[g],\si) = \inf_{\tilde{g} \in [g] } \lambda_1^+(\tilde{g})
  \Vol(M,\tilde{g})^{1/n}   $.
On montre que
$0< \lamin(M,[g],\si) \leq \lamin(\mS^n)$.
On trouve des conditions suffisantes pour lesquelles on obtient
l'in\'egalit\'e stricte    $\lamin(M,[g],\si) < \lamin(\mS^n)$.
Cette in\'egalit\'e stricte a des applications en g\'eom\'etrie
spinorielle conforme.
\end{abstract}
\end{frontmatter}


\noindent{\bf Abridged English Version}\\

\noindent Let $(M,g,\si)$ be a compact spin manifold of dimension $n \geq 2$.
For a metric $\tilde g$ in the conformal class $[g]$ of $g$,
let $\lambda_1^+(\tilde{g})$ be the smallest positive
eigenvalue of the Dirac operator $D$ in the metric $\tilde{g}$. Similarly, let
$\lambda_1^-(\tilde{g})$ be the largest negative
eigenvalue of $D$.  We define
$$\lamin(M,[g],\si) = \inf_{\tilde{g} \in [g] } \lambda_1^+(\tilde{g})
\Vol(M,\tilde{g})^{1/n}
\hbox{ and }
\laminm(M,[g],\si) =  \inf_{\tilde{g} \in [g] } |\lambda_1^-(\tilde{g})|\;
\Vol(M,\tilde{g})^{1/n} $$

\noindent As a first result we obtain
\begin{equation}\label{ineq.large+}
  \lamin(M,[g],\si) \leq \lamin(\mS^n) = \frac{n}{2}\, \om_n^{{1 \over n}}
\end{equation}
and
\begin{equation}\label{ineq.large-}
\laminm(M,[g],\si)\leq \lamin(\mS^n)=\lamin(\mS^n).
\end{equation}

\noindent The main result is then the following:

\begin{theorem} \label{main1}
Inequality (\ref{ineq.large+}) (resp. inequality
(\ref{ineq.large-})) is strict 
if $(M,g)$ is conformally flat, if $D$ is invertible and if the
mass endomorphism possesses a positive (resp. negative) eigenvalue.
\end{theorem}

\noindent This result has applications to conformal spin geometry.

\subsection*{Publication de cet article} Cet article
est publi\'e sous la r\'ef\'erence:\\
B. Ammann, E. Humbert, B. Morel, {\it
Un probl\`eme de type Yamabe sur les vari\'et\'es spinorielles compactes.}
C. R. Math. Acad. Sci. Paris {\bf 338} (2004), 929--934.

Apr\`es publication, nous avons trouv\'e une erreur dans l'un de nos r\'esultats
qui implique que certains \'enonc\'es de la version publi\'ee 
ne sont pas corrects. Dans la version actuelle, ces \'enonc\'es sont 
corrig\'es.

\section{Introduction}
Soit $(M,g,\si)$ une vari\'et\'e spinorielle compacte de dimension $n \geq
2$. Si $\tilde{g} \in [g]$ est une m\'etrique conforme \`a $g$, on note
$\lambda_1^+(\tilde{g})$ (resp. $\lambda_1^-(\tilde{g})$) la plus petite
valeur propre positive (resp. n\'egative) de
l'op\'erateur de Dirac par rapport \`a la  m\'etrique $\tilde{g}$.
On d\'efinit alors
$$\lamin(M,[g],\si) = \inf_{\tilde{g} \in [g] } \lambda_1^+(\tilde{g})
\Vol(M,\tilde{g})^{1/n}
\qquad\hbox{et}\qquad
 \laminm(M,[g],\si) =  \inf_{\tilde{g} \in [g] } |\lambda_1^-(\tilde{g})|\;
\Vol(M,\tilde{g})^{1/n} $$
Lott et le premier auteur ont
montr\'e dans \cite{lott:86,ammann:03} que
$$\lamin(M,[g],\si) >0 \hbox{ et } \laminm(M,[g],\si) >0$$
De nombreux travaux sont consacr\'es \`a l'\'etude de ces invariants
conformes, en particulier \cite{hijazi:86}, \cite{lott:86}, \cite{baer:92b},
\cite{ammann:habil,ammann:p03a}.
Les probl\`emes que nous abordons dans cette note sont issus d'un r\'esultat
de Ammann \cite{ammann:habil,ammann:p03a} qui dit que  si $n \not= 2$ ou $Ker(D)= \{
0 \}$ alors
\begin{equation}\label{ineq.large}
  \lamin(M,[g],\si) \leq \lamin(\mS^n) = \frac{n}{2}\, \om_n^{{1 \over n}}
\hbox{ et } \laminm(M,[g],\si)\leq \lamin(\mS^n)=\laminm(\mS^n),
\end{equation}
o\`u $\om_n$ est le volume de la sph\`ere unit\'e standard $\mS^n$.
\`A la lecture de ce r\'esultat se posent plusieurs questions
naturelles: d'abord, est-ce que les in\'egalit\'es (\ref{ineq.large}) sont vraies
dans le cas $n=2$ et $D$ non inversible? Ensuite, \`a quelles conditions
sont-elles strictes?
L'int\'er\^et de r\'epondre \`a cette derni\`ere question est multiple.
En effet, les in\'egalit\'es strictes
\begin{equation} \label{i+strict}
 \lamin(M,[g],\si) < \lamin(\mS^n)
\end{equation}
et
\begin{equation} \label{i-strict}
 \laminm(M,[g],\si) < \laminm(\mS^n) = \lamin(\mS^n)
\end{equation}
ont plusieurs applications. L'une d'elles, via l'in\'egalit\'e
d'Hijazi (\cite{hijazi:86} et \cite{hijazi:91}), a trait au c\'el\`ebre
probl\`eme de Yamabe dont la r\'esolution fut obtenue par  Aubin
\cite{aubin:76} et Schoen \cite{schoen:84} au milieu des ann\'ees
1980. Pour plus d'informations sur ce probl\`eme, le lecteur pourra se r\'ef\'erer \`a
\cite{aubin:76}, \cite{hebey:97} ou \cite{lee.parker:87}.


Une autre application, plus nouvelle, concerne l'\'equation 
$$D(\phi)= \lambda |\phi|^{{2 \over n-1}} \phi, \hbox{ avec }
{\parallel \phi \parallel}_{{2n \over n-1}}=1 $$
Cette \'equation fut \'etudi\'ee dans  \cite{ammann:habil,ammann:p03a} o\`u
il est d\'emontr\'e que si   (\ref{i+strict}) (resp.
(\ref{i-strict})) est vraie alors on peut r\'esoudre l'\'equation
pr\'ec\'edente avec
$\lambda = \lamin(M,[g],\si)$ (resp. $=\laminm(M,[g],\si)$).
Cette EDP est invariante par un changement conforme de m\'etrique, et elle est
critique du point de vue des injections de Sobolev en ce sens
que les injections de Sobolev qui interviennent dans sa r\'esolution
ne sont pas compactes. De ce point de vue, cette \'equation est tr\`es
proche de celle qui intervient
dans la r\'esolution du probl\`eme de Yamabe. Dans
\cite{ammann:habil,ammann:p03a} on montre aussi
que $\lamin(M,[g],\si)$ (resp. $\laminm(M,[g],\si)$)
est atteinte. Il est \`a noter que le sch\'ema de d\'emonstration de ces
r\'esultats est proche de celui utilis\'e pour la r\'esolution du
probl\`eme de Yamabe m\^eme si de nombreuses difficult\'es interviennent en
travaillant avec les spineurs.

Les r\'esultats que nous obtenons sont r\'esum\'es dans le th\'eor\`eme
suivant:
\begin{theoreme} \label{main}
Soit $(M,g,\si)$ une vari\'et\'e spinorielle compacte de dimension $n
\geq 2$. Alors les in\'egalit\'es larges (\ref{ineq.large}) sont toujours
vraies. De plus,
si $(M,g)$ est conform\'ement plate, si $D$ est inversible
et si l'endomorphisme de masse
  (voir dernier paragraphe) poss\`ede une valeur propre strictement
  positive (resp. n\'egative) alors l'in\'egalit\'e (\ref{i+strict})
  (resp. (\ref{i-strict}) ) est vraie.
\end{theoreme}

\begin{remarque}
On montre que l'op\'erateur de masse a un spectre sym\'etrique
  si $n \not\equiv 3 \hbox{ mod } 4$. Dans ce cas, s'il  est non nul,
  on obtient directement que les deux in\'egalit\'es  (\ref{i+strict}) et
  (\ref{i-strict}) sont vraies.
\end{remarque}

\section{Le principe de d\'emonstration du th\'eor\`eme (\ref{main})}
La  d\'emonstration du th\'eor\`eme est bas\'ee sur la construction d'un
spineur test ad\'equat (voir la proposition 1 ci-dessous pour plus de
d\'etails sur cette affirmation). Cela apporte des difficult\'es
suppl\'ementaires par rapport au cas classique des fonctions. Dans le cas
conform\'ement plat, on retrouve par ailleurs un lien \'etroit avec le
th\'eor\`eme de la masse positive utilis\'e par Schoen dans son \'etude 
\cite{schoen:84} du probl\`eme de Yamabe. On d\'efinit dans notre contexte
une masse, l'endomorphisme de masse, comme le terme constant de la fonction
de Green de l'op\'erateur de Dirac.

\noindent La premi\`ere \'etape consiste \`a formuler le probl\`eme sous forme
variationnelle.

\begin{prop}[\cite{ammann:03}]
Soit $\psi\in\Gamma(\Sigma M)$. On d\'efinit
  $$J(\psi)=\frac{\Big(\int_M|D\psi|^{\frac{2n}{n+1}}v_g\Big)^\frac{n+1}{n}}{
\left|\int_M <D\psi,\psi>v_g\right|}$$
Alors
\begin{eqnarray} \label{funct}
\lamin(M,[g],\si) (\hbox{resp. } \laminm(M,[g],\si) ) =\inf_\psi J(\psi)
\end{eqnarray}
o\`u l'infimum est pris sur l'ensemble des spineurs de classe $C^{\infty}$ tels que
$\left(\int_M  <D\psi,\psi>v_g \right)>0$ (resp $<0$).
\end{prop}
Le probl\`eme se ram\`ene ainsi \`a trouver un spineur $\psi$ pour lequel
$J(\psi) < \lamin(\mS^n)$.

\subsection{L'in\'egalit\'e large (\ref{ineq.large})}
Dans ce cas, le spineur $\psi$ s'obtient de la mani\`ere suivante: \`a partir d'un
spineur de Killing sur la sph\`ere (qui r\'ealise l'infimum dans la
fonctionnelle ci-dessus), on obtient via la projection st\'er\'eographique de p\^ole $N$ et la covariance conforme de l'op\'erateur de Dirac, un spineur $\phi$ sur $\mR^n$ qui satisfait
 $$D\phi = \frac{n}{2} f \phi$$
o\`u $f(x)= \frac{2}{1+|x|^2}$ est telle que les m\'etriques de $\mS^n \setminus \{N\}$ et de $\mR^n$ satisfassent $g_{\mS^n}= f^2 g_{eucl}$.
On choisit un point $p \in M$. 
Quitte \`a remplacer $g$ par une m\'etrique qui lui est conforme, on montre
qu'on peut supposer que
$\Ric_g(p)=0$ et $\Delta_g(\Scal_g)(p)=0$. 
On trivialise
le fibr\'e des spineurs autour de $p$ en utilisant la
trivialisation de Bourguignon-Gauduchon (voir
\cite{bourguignon.gauduchon:92,ammann.grosjean.humbert.morel:p07}). On obtient alors un
diff\'eomorphisme $\tau: \Sigma U \to \Sigma V$
o\`u $U$ est un voisinage
de $0$ dans $\mR^n$ et $V$ un voisinage de $p$ dans $M$ tels que pour tout
$x \in V$, $\tau$ est une isom\'etrie de $\Sigma_{exp_p(x)} U$ sur
$\Sigma_xV$.
Soit maintenant $\eta$ une fonction de cut-off de classe $C^{\infty}$ au
voisinage de $p$ et $0$ sur $M \setminus V$. On pose
$$\psi_{\ep} = \eta \tau\circ\phi\left(\frac{x}{\ep}\right)$$
On montre alors que si le spineur de Killing de la sph\`ere duquel
on est parti est bien choisi,
alors,
quand  $\epsilon \to 0$,
$J(\psi_{\ep})=  \lamin(\mS^n) +o(1)$, et on obtient (\ref{ineq.large}).


\subsection{Le cas conform\'ement plat}
Dans ce paragraphe, on suppose que $(M,g)$ est conform\'ement plate et que $D$ est inversible.
On aura alors besoin d'introduire la notion
d'endomorphisme de masse. Soit $p \in M$ un point de $M$ fix\'e. Quitte \`a
faire un changement conforme de m\'etrique, on peut supposer que $g$ est plate au voisinage de $p$.

\subsubsection{L'endomorphisme de masse}
L'endomorphisme de masse est le terme constant dans la fonction de Green de
$D$. La fonction de Green de $D$ est un endomorphisme
$$G_D: \Sigma_pM \to \Gamma(\Sigma (M \setminus \{p\}))$$
qui \`a un \'el\'ement $\psi_0$ de $\Sigma_pM$ associe un spineur  de classe $C^{\infty}$ d\'efini sur $M \setminus \{p\}$ et qui satisfait
$$D (G_D(\psi_0)) = \psi_0 \delta_p$$
au sens des distributions.
Ici, $\delta_p$ est la masse de Dirac en $p$. Alors,
\begin{prop}\label{expan}\
La fonction de Green $G_D$ existe et est unique. De plus, si on choisit une carte $\phi = (x_1,...,x_n)$ dans laquelle la m\'etrique est euclidienne et telle que $\phi(p)=0$, alors on peut identifier $\Sigma U$ ($U$ voisinage de $0$ dans $\mR^n$) et
$\Sigma V$ ($V$ voisinage de $p$ dans $M$ ) de mani\`ere triviale, gr\^ace par exemple \`a la trivialisation de Bourguignon-Gauduchon.
Avec cette identification, $G_D$ a le d\'eveloppement suivant quand $x$ tend vers $p$
$$\omega_{n-1} \,G_D \psi_0 (x)= - \frac{x}{|x|^n}  \cdot \psi_0 +
v(x)\psi_0$$
o\`u $v(x)(\psi_0)$ est un spineur harmonique d\'efini au voisinage de $p$.
\end{prop}

\begin{definition}
On appelle \emph{endormophisme de masse} l'endomorphisme suivant:
\[ \alpha: \left| \begin{array}{ccc}
\Sigma_p M & \to & \Sigma_p M \\
\psi_0 & \mapsto & v(p)(\psi_0)
\end{array} \right. \]
\end{definition}

\noindent L'endomorphisme de masse est lin\'eaire et
autoadjoint. Ses valeurs propres sont donc r\'eelles.
Si $n \not\equiv 3 \hbox{ mod } 4$, son spectre est sym\'etrique.
De plus, le signe de ses valeurs propres ne d\'epend pas de la m\'etrique
plate au voisinage de $p$ choisie dans la classe conforme de $g$.

\begin{exemple}

\noindent - Si $(M,g)$ est un tore plat alors l'endomorphisme de masse est nul.

\noindent - Si $(M,g)$ est l'espace projectif r\'eel de dimension $n \equiv
3 \hbox{ mod } 4$, alors l'endomorphisme de masse est un multiple non nul
de l'identit\'e.
\end{exemple}

\begin{remarque}
Supposons $n$ impair. Alors, pour tout $z\in\mC,
\Re(z)<-n$, le noyau d'int\'egration $k(z;x,y)\in \Hom(\Si_xM,\Si_yM)$
de $D^z$ est continu en les variables $(x,y)$ sur $M \times M$.   La fonction
$h(z;x):=k(z;x,x)$ a une extension m\'eromorphe pour  $z\in\mC$.
K. Okikiolu \cite[Theorem 1.2.2]{okikiolu:01} d\'emontre que
$h(-1;p)=\alpha$.
\end{remarque}

\subsection{Preuve du Th\'eor\`eme \ref{main} dans le cas conform\'ement plat}

Soit $\psi_0$ un vecteur propre de $\alpha$ associ\'e \`a la valeur propre $\nu$. Soit aussi $\ep>0$ petit et
 $$\rho:=\varepsilon^{\frac{1}{n+1}}\qquad \ep_0:=\frac{\rho^n}{\varepsilon} f(\frac{\rho}{\varepsilon})^{\frac{n}{2}}\;.$$
Les spineurs-tests $\psi_{\ep}^+$ (pour $\lamin$) et $\psi_{\ep}^-$ (pour $\laminm$) utilis\'es ici sont les suivants

\[ \psi_\ep^\pm:= \left|
\begin{array}{ccc}
f(\frac{r}{\varepsilon})^{\frac{n}{2}}
(1\mp\frac{x}{\varepsilon})\cdot\psi_0\mp
\varepsilon_0 \nu \psi_0  & \hbox{ si } &  r\leq\rho\;, \\[5mm]
\mp\varepsilon_0(G_D(\psi_0) -\eta v(x)(\psi_0) - \nu \psi_0)+\eta
f(\frac{\rho}{\varepsilon})^{\frac{n}{2}}\psi_0 & \hbox{ si } & \rho\leq\
r\leq 2\rho\;,\\[5mm]
\varepsilon_0 \,G_D(\psi_0) & \hbox{ si }&  r\geq 2\rho\;,
\end{array}      \right. \]
o\`u  $f$ est comme au paragraphe pr\'ec\'edent, o\`u
$r=|x|$ et  o\`u $\eta$ est une fonction de cut-off function \'egale
\`a $1$ sur $B(p,\rho)$, \'egale \`a $0$ sur le compl\'ementaire de  $B(p,2\rho)$ et qui satisfait
  $|\nabla\eta|\leq\frac{2}{\rho}$. Il faut remarquer que le spineur $\phi$
  construit au paragraphe pr\'ec\'edent exprim\'e en coordonn\'ees dans
  $\mR^n$ a la forme suivante:
$$\phi= f(\frac{r}{\varepsilon})^{\frac{n}{2}}(1\mp\frac{x}{\varepsilon})\cdot\phi_0$$
o\`u $\phi_0$ est un spineur constant. Donc le spineur $\psi_{\ep}$, comme dans le cas non conform\'ement plat est \'egal \`a $\phi(\frac{x}{\ep})$ dans un voisinage de $p$.
On calcule alors que
$$J(\psi_{\epsilon} ) \leq  \lamin(\mS^n) + C_0 \nu \epsilon^{n-1} + o( \epsilon^{n-1})$$
o\`u $C_0$ est une constante
strictement positive.
Pour $\epsilon$ petit, on trouve que $J(\psi_{\ep}) < \lamin(\mS^n)$. Cela
prouve le th\'eor\`eme (\ref{main}).


\def\refname{R\'ef\'erences}







\end{document}